\documentclass[12pt]{article}
\usepackage{url}
\title{\textbf{Algorithmic Information Theory\\Some Recollections}}
\author{Gregory Chaitin, 7 January 2007}
\date{}
\begin{document}
\maketitle
\sloppy

\section*{Introduction}

AIT is a theory that uses the idea of the computer, particularly the size of computer
programs, to study the limits of knowledge, in other words, what we can know, and how.
This theory can be traced back to Leibniz in 1686, and it features a place in pure mathematics
where there is absolutely no structure, none at all, namely the bits of the halting
probability $\Omega$.

There are related bodies of work by other people going in other directions,
but in my case the emphasis is on using the idea of algorithmic complexity to obtain incompleteness results.
I became interested in this as a teenager and have worked on it ever since.

Let me tell you that story.
History is extremely complicated, with many different points of view.
What will make my account simple is the unity of purpose imposed on a field
that is a personal creation, 
that has a central spine, that pulls a single thread.
What did it feel like to do that? In fact, it's not something I did.
It's as if the ideas wanted to be expressed through me.
    
It is an overwhelming experience to feel possessed by promising new ideas.  
This happened to me as a teenager, and I have spent the rest of my life
trying to develop the ideas that flooded my mind then. These ideas were deep enough to merit
45 years of effort, and I feel that more work is still needed.
There are many connections with crucial concepts in other fields: physics, biology, philosophy, theology,
artificial intelligence\ldots\ Let me try to remember what happened to me\ldots\ The history of a person's life,
that's just gossip. But the history of a person's ideas, that is real, that is important, that is
where you can see creativity at work. That is where you can see new ideas springing into being.

\section*{AIT in a Nutshell}

G\"odel discovered incompleteness in 1931 using a version of the liar paradox, ``This statement is unprovable.''
I was fascinated by G\"odel's work.
I devoured Nagel and Newman, \emph{G\"odel's Proof}, when it was published in 1958.

I was also fascinated by computers, and by the computer as a mathematical concept.
In 1936 Turing derived incompleteness from uncomputability.
My work follows in Turing's footsteps, not G\"odel's, but adds the idea
of looking at the size of computer programs.

For example, let's call a program $Q$ ``elegant'' if no program written in the same language
that is smaller than $Q$ produces the same output.
Can we prove that individual programs are elegant?
In general, no.
Any given formal axiomatic system can only enable us to show that finitely many programs are elegant.

It's easy to see that this must be so. Just consider a program $P$ that calculates
the output of the first provably elegant program that is larger than $P$.
$P$ runs through all the possible proofs in the formal axiomatic system until it finds the first
proof that an individual program $Q$ larger than $P$ is elegant, 
and then $P$ runs $Q$ and returns $Q$'s output as its ($P$'s) output.

If you assume that only true theorems can be proved in your formal axiomatic system, then
$P$ is too small to be able to produce the same output as $Q$.
If $P$ actually succeeds in finding the program $Q$, then we have a contradiction.
Therefore $Q$ is never found, which means that no program that is bigger than $P$ can be proven to be elegant.

So how big is $P$?
Well, it
must include a big subroutine for running through all the possible proofs of the formal axiomatic system.
The rest of $P$, the main program, is rather small; $P$ is mostly that big subroutine.
That's the key thing, to focus on the number of bits in that subroutine.

So let's define the algorithmic complexity of a formal axiomatic system to be
the size in bits of the smallest program for running through all the proofs and producing all the theorems.
Then we can state what we just proved like this:
You can't prove that a program is elegant if its size is 
substantially
larger than the algorithmic complexity of the formal axiomatic 
system that you are using.

Instead of saying ``a formal axiomatic system of algorithmic complexity $N$,'' I'll just say ``$N$ bits of axioms.''
So
if you have $N$ bits of axioms, then no program larger than $N+c$ bits in size can be proven to be elegant.
That's the result we just proved.

A more sophisticated example is the number I call $\Omega$, which is the halting probability of a computer
running a program produced one bit at a time by repeatedly tossing a coin.
Because it is a probability,
this number has to be between zero and one.
Imagine writing it out in binary:
\[
   \Omega = .011100\ldots
\]
These bits are peculiar, they are
irreducible mathematical information.
This means that
a formal axiomatic system with $N$ bits of axioms can enable you to determine at most $N+c$ bits of $\Omega$.
Essentially the only way to determine bits of $\Omega$ is to add that information directly to your axioms.
Even though $\Omega$ is a single well-defined real number (once you fix the programming language),
its bits have no structure, no pattern, none at all, they are irredundant, irreducible mathematical information.

In other words,
the bits of $\Omega$ are mathematical facts that are true for no reason, no reason simpler than themselves.

So that's the basic idea, and those are my two favorite results, but the devil is in the details.
You can spend your life on those details, and I did.

\section*{Chaitin Research Timeline}

\begin{itemize}

\item\textbf{1947}: 
Born in Chicago, child of Argentine immigrants. Family moves to New York.
      
\item\textbf{1956}: 
Nagel and Newman's article on ``G\"odel's proof'' is published in \emph{Scientific American}.
Article contains a photo by Arnold Newman of G\"odel sitting in front of an empty
blackboard at the Princeton Institute for Advanced Study.
   
\item\textbf{1958}: 
Nagel and Newman's book \emph{G\"odel's Proof}
is published by New York University Press.

\item\textbf{1959}: 
Following directions in the \emph{Scientific American} ``Amateur Scientist'' department, I build a Van de Graaff
generator for high-voltage static electricity. 

\item\textbf{1962}: 
First year at Bronx High School of Science.
While answering an essay question on the entrance exam for the Columbia University Science Honors Program 
for bright high school students,
I get the idea of defining randomness using program-size complexity. 

The essay question is what do you conclude if you find a pin on the moon.\footnote
{This was before the first lunar landing.}
My answer is that this means that somebody must have visited before you, because a pin is not
natural, it is artificial, the product of intelligence.  And, I remark, what this means is
that there is a small program to calculate it, to create it. That's how we can tell that the pin
has structure and is artificial. And, contrariwise, something
natural would not have a description that can be compressed into a small program, because it was not
designed.

And then, as a throw-away remark, I state that a random thing is one that cannot 
be compressed into a smaller program. More precisely, I am speaking about a digital description of an
object, not about the object itself. 
In other words, in 1962 I give the following
\begin{itemize}
\item
\textbf{Definition of Randomness R1}:
A random finite binary string is one that cannot be
compressed into a program smaller than itself,
that is, that is not the unique output of a program without any input, a program 
whose size in bits is smaller than the size in bits of its output.
\end{itemize}
However I quickly
forget about this definition, because I am having so much fun
learning how to write, debug and run computer programs in the Science Honors Program.
And I am given the run of the math stacks at Columbia University and can
hold in my hands and study the collected works of Euler and other priceless volumes.

\item\textbf{1963}:
Shannon and McCarthy, \emph{Automata Studies}, Princeton University Press, 1956,
contains E. F. Moore's paper ``Gedanken-experiments on sequential machines.''\footnote
{I became aware of Shannon and McCarthy,
perhaps the first book on the theory of computation,
because it was reviewed in \emph{Scientific American}.}
Following Moore, I write a program for
identifying a finite-state black box by putting in inputs and looking at the outputs.
My experiments suggest this is easier to do than Moore anticipated. I prove
this to be the case in a note ``An improvement on a theorem of E. F. Moore''
(\emph{IEEE Transactions on Electronic Computers}, 1965), my first publication.

\item\textbf{1964}:
Summer vacation between high school and college, I 
try to find an infinite set with no subset that is easier to generate than the entire set.
By easier I mean faster or simpler; at this point I am simultaneously exploring run-time complexity and program-size complexity.
The work goes well, but is not
published until 1969 in the \emph{ACM Journal} as ``On the simplicity and speed of programs for
computing infinite sets of natural numbers.''
    
Also that summer, I get the first incompleteness result that I will publish,
\textbf{UB1}, an upper bound on the provable lower bounds on run-time complexity 
in any given formal axiomatic system.
This is published in 1970 in a Rio de Janiero
Pontif\'{\i}cia Universidade Cat\'olica research report, 
and only there.\footnote{While writing up that report in Rio, I realize I can also obtain
an upper bound on the provable lower bounds on program-size complexity.}
     
Another discovery that summer, \textbf{UB2}, is that one can
diagonalize over the output of all programs that provably calculate 
total functions $f\!\!:N\rightarrow N$
to obtain a faster growing computable total function $F\!\!:N\rightarrow N$. 
That is to say, given any formal axiomatic system, one can construct a computer program from it
that calculates a total function $f\!\!:N\rightarrow N$, 
but the fact that this program calculates a total function
$f\!\!:N\rightarrow N$ cannot be proved within the formal axiomatic system, 
because $f$ goes to infinity too quickly.
``Calculates a total function $f\!\!:N\rightarrow N$'' merely means 
that every time we give the program $f$ a natural
number $n$ as input, it eventually outputs a single natural number $f(n)$ and then halts.
     
The result UB1 is actually a corollary of UB2, 
since all lower bounds on run-time complexity are computable total functions.

Now one would say that the proof of UB2 is an instance of Cantor diagonalization, but
in my opinion
it's really closer to Paul du Bois-Reymond's theorem on orders of infinity. His theorem is that for any
scale of rates of growth, any infinite list of functions that go to infinity faster and faster, 
for example 
\begin{eqnarray*}
   f_0(n) & = & 2^n, \\
   f_1(n) & = & 2^{2^n}, \\
   f_2(n) & = & 2^{2^{2^n}} \ldots,
\end{eqnarray*}
there is another function 
\[
f_{\omega}(n) = \max_{k \le n} f_k(n)
\]
that goes to infinity even more quickly. 
As far as I know, Paul du Bois-Reymond's work was independent of Cantor's. 
     
Note the Cantor ordinal number $\omega$ as a subscript.  We can then form 
\begin{eqnarray*}
f_{\omega+1}(n) & = & 2^{f_{\omega}(n)}, \\
f_{\omega+2}(n) & = & 2^{f_{\omega+1}(n)}, \\
f_{\omega+3}(n) & = & 2^{f_{\omega+2}(n)} \ldots
\end{eqnarray*}
and then
\[
f_{2\omega}(n) = \max_{k \le n} f_{\omega+k}(n).
\]
Continuing in this manner, we get to 
\[
f_{3\omega} \ldots
f_{4\omega} \ldots
f_{\omega^2} \ldots
f_{\omega^3} \ldots
f_{\omega^{\omega}} \ldots
f_{\omega^{\omega^{\omega}}} \ldots
\]
and onwards and upwards, incredibly fast-growing functions all.
     
I had read about Paul du Bois-Reymond's work
in a monograph by G. H. Hardy called \emph{Orders of Infinity}.\footnote{I learned
the calculus from 
Hardy's \emph{A Course of Pure Mathematics}, and also enjoyed
\emph{A Mathematician's Apology} and Hardy and Wright, 
\emph{An Introduction to the Theory of Numbers}.}

My point of view changes
between 1964 and 1965. In my 1964 work,
only \textbf{infinite} computations are considered. In 1965, on the contrary,
only \textbf{finite} computations are considered, computations that produce \textbf{a single output}.
Furthermore, my interest shifts from run-time complexity to program-size complexity.
                                                                                                                                     
\item\textbf{1965}: 
During the spring term of my first year at City College, CUNY, I simultaneously study three books:
von Neumann and Morgenstern, \emph{Theory of Games and Economic Behavior},
Shannon and Weaver, \emph{The Mathematical Theory of Communication}, and
Turing's 1936 paper ``On computable numbers\ldots''\ in the anthology Davis, \emph{The Undecidable}.
The 1962 definition of randomness (R1) comes back to me; as I will now explain, all three books play a vital role.
It all comes together
as I am reading a discussion in von Neumann and Morgenstern of
the game of matching pennies, for which
their theory says that you should toss a coin.\footnote
{One of the players is trying to match the other player's choice of head or tails.}

In a footnote they remark that
the theory of games seems to require a quantum-mechanical world in which God plays dice.
Not really, I say to myself.
Another logical possibility would be that the theory of games tells you to use a random sequence
of choices, but you cannot compute this sequence of choices from the theory.

You see, in the game of matching pennies, if a theory can tell you exactly what to do, you can
predict what your opponent will do and beat him. The solution to the paradox is either
that the theory asks you to use
physical randomness, which is
unpredictable, or that you have a theory that says
you must make mathematically random or unstructured choices, and
the contradiction is avoided because these are in fact uncomputable (a notion taken from Turing).
 
The third leg of the stool comes from reading Shannon,
who defines a message to be random
or have maximum entropy if it cannot be compressed,
if it cannot be encoded more compactly.
Obviously the most general possible \textbf{decoder} would be a universal Turing machine, a general-purpose
computer.
 
With my 1962 definition (R1), I have managed
to connect game theory, information theory and computability theory.
Now all I have to do is work out the details.

Now it's
the summer vacation between my first and second year at City College, 
and I attempt to carry out the plan.
At the beginning of the summer, the road forward seems blocked, but I keep trying.
Later in the summer, ideas begin to flood into my mind.
I write a single paper that is the size of a small book.
At the request of the editors, I later divide it in two,
and delete much material to save space.
     
This paper ``On the length of programs for computing finite binary sequences''---part one is 
published in 1966 and part two is published in 1969, both in the \emph{ACM Journal}---presents 
\textbf{three} different theories of program-size complexity
(and embryonic versions of ideas that I would explore for years):
\begin{itemize}
\item
\textbf{Complexity Theory (A)}: 
Counting the number of states in a normal Turing machine with a fixed number of tape symbols.
I call this Turing machine state-complexity.
\item
\textbf{Complexity Theory (B)}:
The same as theory (A), but now there's
a fixed upper bound on the size of transfers---jumps, branches---between states. 
You can only jump nearby.
I call this bounded-transfer Turing machine state-complexity.
\item
\textbf{Complexity Theory (C)}:
Counting the number of bits in a binary program, a bit string.
The program starts with a \textbf{self-delimiting} prefix, indicating which computing machine to simulate,
followed by the binary program for that machine.
That's how we get what's called \textbf{a universal machine}.\footnote
{I call theory (C) ``blank-endmarker'' program-size complexity, to distinguish it from
``self-delimiting'' program-size complexity, theory (D) below.} 
\end{itemize}

Let's define the complexity of a bit string to be the size of the smallest program that computes it.
     
In each case, theory (A), (B) or (C), I show that most $n$-bit strings have complexity close to the maximum possible, 
and I determine asymptotic formulas for the maximum possible complexity of an $n$-bit string. 
These maximum or near maximum complexity strings are defined to be random. 
To show that this is reasonable, I prove, for example, that these strings are ``normal'' in Borel's sense. 
This means that
all possible blocks of 
bits of the same size occur in such strings approximately the same number of times,
an equi-distribution property.
    
I start with theory (A) because that seems the most straightforward thing to do.
The idea in theory (A) is to eliminate all the redundancy in a real programming language.
Then I switch to theory (B), in which I don't eliminate the redundancy, I live with it.
The proofs are prettier; more subtle, not so heavy-handed.
However, in theories (A) and (B) I cannot figure out how to show that
a \textbf{small} amount of structure in an $n$-bit string will force
its complexity to dip below the maximum possible complexity for $n$-bit strings.

To solve this,
I switch from Turing machines to binary programs and theory (C).
Theory (C) solves the problem, but feels too easy, like stealing candy from a baby.
For example, in theory (C) it is trivial to show that most $n$-bit strings have close to the maximum possible complexity.
And this maximum possible complexity is precisely $n+1$, not an asymptotic estimate as in theories (A) and (B).\footnote
{To get (max complexity $n$-bit string) $= n + 1$, theory (C) has to be a bit more complicated.
\begin{itemize}
\item
\textbf{Complexity Theory (C2)}: If the first bit of the program is a 0, then output the rest of the program as is and halt.
If the first bit of the program is a 1, process the rest of the program as in theory (C).
\end{itemize}
(C2) is the version of theory (C) given in ``On the length of programs for computing finite binary sequences:
statistical considerations'' (\emph{ACM Journal}, 1969), the second of the two papers put together from my
1965 randomness manuscript.}
     
By the way, in theories (A) and (B), randomness definition (R1) does not apply, because the size of programs is measured
in states, not bits.  It is necessary to use a slightly different definition of randomness:
\begin{itemize}
\item
\textbf{Definition of Randomness R2}:
A random $n$-bit string is one
that has maximum or near maximum complexity.
In other words, 
an $n$-bit string is random
if its complexity is approximately equal to the maximum complexity of any $n$-bit string.
\end{itemize}
In theory (C), (R1) works fine (but so does (R2), which is more general).
      
Theory (C) is essentially the same as the one independently proposed by Kolmogorov at
about the same time (1965).\footnote
{Solomonoff was the first person to publish the idea of program-size complexity---in fact, 
(C)---but he did not propose a definition of randomness.}

However, I am dissatisfied with theory (C); the absence of \textbf{subadditivity} disturbs me.
What is subadditivity? 
The usual definition is that a function $f$ is subadditive if $f(x + y) \leq f(x) + f(y)$. I mean something slightly different.
Subadditivity holds if
the complexity of computing two objects
together (also known as their \textbf{joint complexity}) is bounded by the sum of their individual complexities.\footnote
{In the case of joint complexity the computer has \textbf{two} outputs, or outputs \textbf{a pair} of objects, 
whatever you prefer.}
In other words, subadditivity means that you can combine subroutines by concatenating them,
without having to add information to indicate where the first subroutine ends and the second one begins.
This makes it easy to construct big programs.
Complexity is subadditive in theories (A) and (B), but not in theory (C).

Last but not least, ``On the length of programs for computing finite binary sequences'' contains 
what I would now call a Berry paradox proof
that program-size complexity is uncomputable. This seed was to grow into my 1970 work on incompleteness, where I
refer to the Berry paradox explicitly for the first time.

\item\textbf{1966}: 
Awarded by City College the Belden Mathematical Prize and the Gitelson Medal ``for the pursuit of truth.''
Family moves back to Buenos Aires.  

\item\textbf{1967}: 
I join IBM Argentina, working as a computer programmer.

\item\textbf{1969}: 
Stimulated by von Neumann's posthumous \emph{Theory of Self-Reproducing Automata},
I work on a mathematical definition of life 
using program-size complexity. This
is published in Spanish in Buenos Aires, 
and the next year (1970) in English in the \emph{ACM SICACT News}. 
This is the first of what on the whole I regard as an unsuccessful series of papers on
theoretical biology.\footnote
{The latest one is ``Speculations on biology, information and complexity''
(CDMTCS Research Report No.\ 282, Univ.\ of Auckland, July 2006).}

\item\textbf{1970}: 
I visit Brazil and inspired by this tropical land, I realize that one can get
powerful incompleteness results using program-size arguments. 
In fact, one can place \textbf{upper} bounds on the provable \textbf{lower} bounds
on run-time \textbf{and} program-size complexity in a formal axiomatic system.
And this provides a way to measure the power of that formal axiomatic system.

This
first information-theoretic incompleteness result is immediately published in a Rio de Janiero 
Pontif\'{\i}cia Universidade Cat\'olica 
research report and also as an \emph{AMS Notices} abstract, and comes out the next year (1971) as
a note in the \emph{ACM SIGACT News}. 

I obtain a \emph{LISP 1.5 Programmer's Manual} in Brazil and start
writing LISP interpreters and inventing LISP dialects.\footnote
{Around 1973, I give courses on LISP and on computability and metamathematics
at the University of Buenos Aires.}

\item\textbf{1971}: 
I write a longer paper on incompleteness,
``Information-theoretic limitations of formal systems,''
which is 
presented at the Courant Institute Computational Complexity Symposium in New York City
in October 1971.
A key idea in this paper is to measure the complexity of a formal axiomatic
system by the size in bits of the program that generates all of the theorems by systematically running
through the tree of all possible proofs. 

\item\textbf{1973}: 
I complete a greatly expanded version of 
``Information-theoretic limitations of formal systems.''
The expanded version appears
in the \emph{ACM Journal} in 1974.
A less technical paper on the same subject, ``Information-theoretic computational complexity,''
is presented at 
the IEEE International Symposium on Information Theory, in Ashkelon, Israel, June 1973,
and is
published in 1974 as an invited paper in the \emph{IEEE Transactions on Information Theory.}\footnote
{In 1974 I send
a copy of this paper to Kurt G\"odel, leading to a pleasant, short phone conversation with
G\"odel and an appointment to meet him at the Princeton Institute for Advanced Study,
an appointment that G\"odel cancels at the last minute.}

\item\textbf{1974}: 
I am invited to visit the IBM Watson Lab in Yorktown Heights for a few months.
The visit goes well, with a number of major breakthroughs.
I realize what to do to theory (C) to restore subadditivity, and 
discover the halting probability $\Omega$.
\begin{itemize}
\item
\textbf{Complexity Theory (D)}:
Counting the number of bits in a self-delimiting binary program, a bit string with the property that
you can tell where it ends by reading it bit by bit without ever reading a blank endmarker.
Now a program starts with a self-delimiting prefix as before, but the program to be simulated
that follows the prefix must \textbf{also} be self-delimiting.
So the
idea is that the \textbf{whole} program must now have the same property the \textbf{prefix} already had in theory (C).
\end{itemize}
                                                                                                                                     
(D) is the mature theory, I believe.
I immediately put this out as an IBM research report.
I present this paper at the opening plenary session of
the IEEE International Symposium on Information Theory in Notre Dame, Indiana, in 
October 1974.
It is published in the \emph{ACM Journal} in 1975
as ``A theory of program size formally identical to information theory.''   
     
There are three key ideas in this paper: self-delimiting programs, a new definition of relative complexity,
and the idea of getting program-size results \textbf{indirectly}
from probabilistic, measure-theoretic arguments
involving the probability $P(x)$ that a program will calculate $x$. 
I call this the algorithmic probability of $x$.\footnote
{Solomonoff tried to define $P(x)$ but could not get $P(x)$ to converge 
since he wasn't working with self-delimiting programs.}
Summing $P(x)$ over all possible outputs $x$ yields the halting probability $\Omega$:
\[
      \Omega = \sum_x P(x).\footnote
{This definition is a bit abstract. Here are two other ways of defining an $\Omega$ number.
As a sum over programs $p$:
\[
      \Omega' = \sum_{\mbox{\scriptsize $p$ halts}} 2^{-|p|}.
\]
As a sum over all positive integers $n$:
\[
      \Omega'' = \sum_n 2^{-H(n)}.
\]
Here $|p|$ is the size in bits of the program $p$, and $H(n)$ is the size in bits of the smallest program for
calculating the positive integer $n$.}
\]
And a key theorem 
\[
\mathbf{(*)} \;\;\; H(x) = -\log_2 P(x) + O(1)
\]
permits us to translate complexities into probabilities and vice versa.
Here the complexity $H(x)$ of $x$ is the size in bits of the smallest program for calculating $x$,
and the $O(1)$ indicates that the difference between the two sides of the equation is bounded.

Incidentally, $(*)$ implies 
that there are few minimum or near-minimum size
programs for calculating something, few minimal descriptions.
That is, an elegant program for calculating something is essentially unique.\footnote
{For more on this, see my book \emph{Exploring Randomness} (2001).
I had proven this result in theory (C) in 1972, 
but that proof wasn't published until 1976, in 
``Information-theoretic characterizations of recursive infinite strings'' 
in \emph{Theoretical Computer Science}.}
    
Where did the halting probability
$\Omega$ come from? How did I come up with it?  Well, already in part two of my first paper on randomness,
``On the length of programs for computing finite binary sequences:
statistical considerations'' (\emph{ACM Journal}, 1969), I use a Heine-Borel-style  
algorithm to exhibit a specific example of a random
infinite sequence of bits.  I think it is important to come
up with specific examples.

The halting probability $\Omega$ is a natural example of a random infinite sequence of bits. 
Besides providing a connection with the work of Turing, 
$\Omega$ makes randomness more concrete and more believable.

Furthermore, once you have a natural example of randomness, you immediately get an incompleteness theorem
from it, as I point out in ``G\"odel's theorem and information'' (\emph{International Journal of Theoretical Physics}, 1982).
My work in 1987 on $\Omega$ and its diophantine equation makes this fully explicit. 
I had been aware of this opportunity for getting a dramatic incompleteness result for some time.\footnote
{See the end of the introductory section of  
``Information-theoretic limitations of formal systems'' (\emph{ACM Journal}, 1974).}
    
$(*)$ is nice but it is only part of the story.
The true reward for changing from theory (C) to (D)
is this spectacular result:
\[
\mathbf{(**)} \;\;\; H(x,y) = H(x) + H(y|x) + O(1).
\]
In words, (the joint complexity of two objects) is equal to the sum of (the absolute complexity of the
first object) plus (the relative complexity of the second object given the first object).
    
To get $(**)$, self-delimiting programs aren't enough, you also need
the right definition of \textbf{relative complexity}.\footnote
{Levin claims to have published theory (D) first. 
However he missed this vital part of theory (D).}
I had used relative complexity in
my big 1965 randomness manuscript, but had eliminated it to save space.
In my 1975 \emph{ACM Journal} paper I take up relative complexity again, but define
$H(x|y)$, the complexity of $x$ given $y$, 
to be the size in bits of the smallest self-delimiting program for calculating
$x$ if we are given for free, not $y$ directly, but a minimum-size self-delimiting program for $y$.
    
And $(**)$ implies that the extent to which computing two things together
is cheaper than computing them separately, also known as the mutual information 
\[
H(x:y) \equiv H(x) + H(y) - H(x,y),
\]
is essentially the same, within $O(1)$, of the extent to which knowing $x$ helps us to know $y$ 
\[
H(y) - H(y|x),
\]
and this in turn is essentially the same, within $O(1)$, 
of the extent to which knowing $y$ helps us to know $x$
\[
H(x) - H(x|y).
\]
This is so pretty that I decide never to use theory (C) again.
For (D) doesn't just restore subadditivity to (C), it reveals an elegant new landscape with sharp results
instead of messy error terms.
From now on, theory (D) only.

\item\textbf{1975}: 
My first \emph{Scientific American} article, ``Randomness and mathematical proof,''
appears. I move back to New York and join the IBM Watson Lab.

In the paper ``Algorithmic entropy of sets'' (\emph{Computers \& Mathematics with Applications}, 1976,
written at the end of 1975), 
I attempt to extend the self-delimiting approach to programs for generating
infinite sets of output.  Much remains to be done.\footnote
{If this interests you,
please see the discussion of infinite computations in the last chapter of \emph{Exploring Randomness} (2001).}

This topic is important, because I think of a formal axiomatic system as a computation
that produces theorems.  My measure of the complexity of a formal axiomatic system
is therefore the size in bits of the smallest \textbf{self-delimiting} program for generating the infinite set of theorems.

\item\textbf{1976--1985}: 
I concentrate on software and hardware engineering for IBM's
RISC (Reduced Instruction Set Computer) project.
    
Even though I spend most of my time on this engineering project,
in 1982 I publish ``G\"odel's theorem and information'' in the \emph{International Journal of
Theoretical Physics}.\footnote 
{At roughly the same time I fulfill a childhood dream by building my own telescope:
I join a telescope-making club and
grind a 6 inch f/8 mirror for a Newtonian reflector
in a basement workshop at the Hayden Planetarium of
the Museum of Natural History.}
This is later included with my 1974 
\emph{IEEE Transactions on Information Theory} 
paper in the influential anthology Tymoczko, \emph{New Directions in the Philosophy of Mathematics},
Princeton University Press, 1998.

My 1985 publication ``An APL2 gallery of mathematical physics: a course outline'' gives computational
working models of physical phenomena to be inserted between chapters of Einstein and Infeld, 
\emph{The Evolution of Physics}.
This APL2 physics simulation software is extremely concise.\footnote
{Besides learning the physics and some numerical analysis, 
I wanted to get a feel for the algorithmic complexity of the laws of physics.}

\item\textbf{1986}: 
My RISC engineering work stops because of an invitation by Cambridge University Press 
to write the first volume in their series Cambridge Tracts in Theoretical
Computer Science.  
I start working on the book, intending merely to collect previous results,
but then the flow of new ideas resumes.
      
Some of these new ideas are presented in the paper
``Incompleteness theorems for random reals'' (1987). 
This paper contains a proof of the
basic result that an $N$-bit formal axiomatic system cannot enable you to determine more than $N + c$ bits
of $\Omega$, bits that may be scattered and do not have to be together
at the beginning of $\Omega$. This is the measure-theoretic proof that I give
in the Cambridge book.
After writing this paper, work on the book begins in earnest.
    
Using work by Jones and Matijasevic on Hilbert's 10th problem,
I calculate a 200-page diophantine equation for $\Omega$.\footnote
{It's actually what's called an \textbf{exponential} diophantine equation.} 
This equation has thousands of unknowns
and a parameter $k$, and has finitely or infinitely many whole-number solutions depending
on whether the $k$th bit of $\Omega$ is, respectively, a 0 or a 1.
To get this equation,
I convert a register machine program for a LISP interpreter into a diophantine equation,
and then I plug into that equation a LISP program for computing lower bounds on $\Omega$.

\item\textbf{1987}: 
Cambridge University Press publishes \emph{Algorithmic Information Theory},
which explains how to obtain the diophantine equation for $\Omega$.
This book also contains a result about random infinite binary sequences. I show that four definitions of
this concept are equivalent: two constructive measure-theoretic definitions due to Martin-L\"of and to
Solovay, and two definitions of my own using program-size complexity.
\emph{Algorithmic Information Theory} is my first publication in which LISP appears.

Simultaneously,
World Scientific publishes a collection of my papers, \emph{Information, Randomness and Incompleteness}.

\item\textbf{1988}: 
I write about the diophantine equation for $\Omega$   
in my second \emph{Scientific American} article, 
``Randomness in arithmetic'' (1988), 
and then in \emph{New Scientist} (1990), and
after that in \emph{La Recherche} (1991).

\item\textbf{1991}: 
I give a lecture on ``Randomness in arithmetic'' in the room where G\"odel taught at the University of Vienna. 

\item\textbf{1992}: 
In the 1992 paper on ``Information-theoretic incompleteness,''
I publish a program-size proof of the theorem that you 
cannot determine the bits of $\Omega$. 
More precisely, an $N$-bit theory---a formal axiomatic system 
with complexity $N$---can permit you to determine at most $N+c$ bits of $\Omega$.
This program-size proof is better than the measure-theoretic proof in the 1987 Cambridge book,
and is the proof that I use in the 1998 book, \emph{The Limits of Mathematics}.

I also publish four papers on LISP program-size complexity: 
\begin{itemize}
\item
\textbf{Complexity Theory (L)}: Counting the number of characters a LISP 
S-expression (that's the program) must have, 
to have a determined value (that's the output).
\end{itemize}
In these four papers several different dialects of LISP are studied,
and appropriate versions of the halting probability 
$\Omega_{\mbox{\scriptsize LISP}}$ are invented for each of them,
together with the corresponding incompleteness theorems.
The techniques
developed in my complexity theories (A) and (B) are put to good use here.

These LISP $\Omega$ numbers may not be as random as the fully random $\Omega$ number in theory (D),
but, like the so-called random infinite sequences in my original 1965 randomness paper,
they come close. For example, they are Borel normal for blocks of all sizes in any base.
    
These 1992 papers are immediately
included in my second World Scientific volume, \emph{Information-Theoretic Incompleteness} (1992). 

I lecture at a meeting on reductionism at Cambridge University. 
The transcript of that lecture, ``Randomness in arithmetic and the decline and fall of reductionism in pure
mathematics,'' appears later in Cornwell, \emph{Nature's Imagination},  Oxford University Press, 1995.

\item\textbf{1995}: 
I discover how to convert theory (D) into a theory about the size of real programs, programs that you can actually run on
a computer (universal Turing machine) that I simulate using a special version of LISP.
\begin{itemize}
\item
\textbf{Complexity Theory (E)}:
Counting the number of bits in a self-delimiting binary program, a bit string with the property that
you can tell where it ends by reading it bit by bit without reading a blank endmarker.
The program starts with a self-delimiting prefix, indicating which computing machine to simulate,
followed by the self-delimiting binary program for that machine, as in theory (D).
But in theory (E), (the prefix 
indicating the machine to simulate)
is a LISP S-expression, 
a program written in a high-level functional programming language, 
that's converted to a bit string.
(E) isn't a new theory, it's a special case of (D)
selected because the prefix indicating the machine to simulate is
encoded in a particularly convenient manner.
\end{itemize}
Now, 30 years after starting to work on program-size complexity, I can finally run programs
and measure their size.
      
Several years are needed to complete this concrete version of AIT and to present it in
my three Springer-Verlag volumes, \emph{The Limits of Mathematics} (1998), \emph{The Unknowable} (1999),
and \emph{Exploring Randomness} (2001).
These books come with LISP software and a LISP interpreter.

Honorary doctorate,  University of Maine.

\item\textbf{2000}: 
Since this year, visiting professor, Computer Science Department, University of Auckland, New Zealand.

\item\textbf{2002}: 
Honorary professor, University of Buenos Aires.

Springer-Verlag publishes 
\emph{Conversations with a Mathematician},
a collection of some of my lecture transcripts and interviews.
       
I'm invited to present a paper at a philosophy congress in Bonn, Germany, in September. 
For this purpose,
I begin to study philosophy, particularly
Leibniz's work on complexity, which I am led to by a hint in a book by Hermann Weyl. 

My paper appears two years later (2004)
in a proceedings volume published by the Academy Press 
of the Berlin Academy that was
founded by Leibniz.
It is reprinted
as the second appendix in my book \emph{Meta Math!}\ (2005).

\item\textbf{2003}: 
Lecture notes \emph{From Philosophy to Program Size} published in Estonia,
based on a course I give there, winter 2003.

\item\textbf{2004}: 
Corresponding member, \emph{Acad\'emie Internationale de Philosophie des Sciences}.

Write \emph{Meta Math!}, a high-level popularization of AIT, published the 
following year by Pantheon Books (2005).
This book is not just an explanation of previous work; it presents a \emph{syst\`eme du monde}.

\emph{Meta Math!}\ 
follows \'Emile Borel in questioning the existence of the bulk of the real numbers,
a train of thought further developed  in my paper ``How real are real numbers?''\
(\emph{International Journal of Bifurcation and Chaos}, 2006).\footnote
{Vladimir Tasi\'c pointed out to me that Borel has a know-it-all real number---a 1927 version of the $\Omega$ number.
My paper on the ontological status of real numbers is dedicated to Borel's memory.}

\item\textbf{2005}: 
I summarize the \emph{syst\`eme du monde} of \emph{Meta Math!}\ in the paper  
``Epistemology as information theory: From Leibniz to $\Omega$'' 
(\emph{\textbf{Collapse}: Journal of Philosophical Research and Development}, 2006).

Honorary president of the scientific committee of the 
Valpara\'{\i}so Complex Systems Institute in Chile.
Member of the scientific advisory panel of FQXi, devoted to Foundational Questions in Physics \& Cosmology.

\item\textbf{2006}: 
The centenary of G\"odel's birth. I publish my third \emph{Scientific American}
article, on ``The limits of reason,'' celebrating Leibniz, whom G\"odel also admired. 
This article is translated and published in about a dozen other languages. 
I summarize my thoughts on incompleteness in an Enriques lecture at the
University of Milan, ``The halting probability $\Omega$: Irreducible complexity in pure mathematics''
(\emph{Milan Journal of Mathematics,} 2007).

A collection of some of my philosophy papers is published in Turin.
Full member, \emph{Acad\'emie Internationale de Philosophie des Sciences}.

\item\textbf{2007}: 
60th birthday.
To paraphrase Einstein, this timeline will have fulfilled its purpose if it shows
how the efforts of a lifetime hang together, and why they lead to certain definite expectations.\footnote
{See the final sentence of Einstein's \emph{Autobiographical Notes}.}
In particular, I think it would be fruitful to explore the following topics.

\end{itemize}

\section*{Challenges for the Future}

\begin{itemize}

\item
On the technical side, many questions remain regarding the program-size complexity
and the algorithmic probability of computing infinite sets of objects.

\end{itemize}

More difficult challenges:

\begin{itemize}

\item
To develop a model of mathematics that is biological, that is, that evolves and develops,
that's dynamic, not static.
Perhaps a time-dependent formal axiomatic system?

\item
To understand creativity in mathematics---where do new ideas come from?---and also in biology---how
do new, much more complicated organisms develop?
Perhaps a life-as-evolving-software model has some merit?

\item
$\Omega=$ concentrated creativity? Each bit of $\Omega=$ one bit of creativity?
Can human intellectual progress be measured by the number of bits of $\Omega$ that we know,
or are currently capable of knowing, as a function of time?

\item
Is the universe discrete or continuous? Can physical systems contain an infinite or only a finite
number of bits?

\item
Make a model world in which you can prove life must develop with high probability.
It doesn't matter if this world isn't exactly like ours; how can that be important?

\item
Is the world built out of information, not matter?
Is it built out of thought?
Is matter just an epiphenomenon, that is, secondary, not primary?
And what are thoughts?

\end{itemize}

\section*{Selected Publications by Chaitin}

\subsection*{Non-Technical Books}
                                                                                                               
\begin{itemize}

\item
\emph{Meta Math!}, Pantheon Books, New York, 2005 (hardcover); Vintage, New York, 2006 (paperback).

\item
U.K. edition: \emph{Meta Maths}, Atlantic Books, London, 2006.

\item
\emph{Conversations with a Mathematician}, Springer-Verlag, London, 2002.\footnote
{Also published in Japanese and Portuguese.}

\end{itemize}

\subsection*{Technical Books}

\subsubsection*{Lecture Notes}
 
\begin{itemize}

\item
\emph{From Philosophy to Program Size}, Institute of Cybernetics, Tallinn, 2003.

\end{itemize}

\subsubsection*{Monographs}

\begin{itemize}

\item
\emph{Algorithmic Information Theory}, Cambridge University Press, 1987 (hardcover), 2004 (paperback).

\item
\emph{The Limits of Mathematics}, Springer-Verlag, Singapore, 1998; 
reprinted by Springer-Verlag, London, 2002.\footnote
{Also published in Japanese.}

\item
\emph{The Unknowable}, Springer-Verlag, Singapore, 1999.\footnote
{Also published in Japanese.}

\item
\emph{Exploring Randomness}, Springer-Verlag, London, 2001.

\end{itemize}

\subsubsection*{Collections of Papers}

\begin{itemize}

\item
\emph{Information, Randomness and Incompleteness}, World Scientific, Singapore, 1987, 2nd ed., Singapore, 1990.

\item
\emph{Information-Theoretic Incompleteness}, World Scientific, Singapore, 1992.

\item
\emph{Teoria algoritmica della complessit\`a}, Giappichelli Editore, Turin, 2006.

\item
\emph{Is God a Computer Programmer?}, e-book, 2007, \\
\url{http://www.cs.auckland.ac.nz/~chaitin/dp2.html}.

\end{itemize}

\end{document}